\documentclass[12pt]{article}

\usepackage{amsmath,amssymb,amsthm}
\usepackage{graphicx}
\usepackage[colorlinks=true,citecolor=black,linkcolor=black,urlcolor=blue]{hyperref}

\usepackage[T1]{fontenc}
\usepackage[latin9]{inputenc}

\usepackage{subfig}
\usepackage{accents}

\usepackage{endnotes}

\usepackage{pstricks,tikz}
\usetikzlibrary{positioning,chains,fit,shapes,calc, fit}
\usetikzlibrary{decorations.pathreplacing,decorations.markings}

\newcommand{\RN}{{\rm RN}}

\theoremstyle{plain}
\newtheorem{theorem}{\protect\theoremname}[section]
  \theoremstyle{plain}
  \newtheorem{proposition}[theorem]{\protect\propositionname}
  \theoremstyle{plain}
  \newtheorem{corollary}[theorem]{\protect\corollaryname}
  \theoremstyle{plain}
  \newtheorem{lemma}[theorem]{\protect\lemmaname}
  \theoremstyle{definition}
  \newtheorem{definition}[theorem]{\protect\definitionname}
 \theoremstyle{plain}
  
 \theoremstyle{plain}

	\makeatother

  \providecommand{\corollaryname}{Corollary}
  \providecommand{\definitionname}{Definition}
  \providecommand{\lemmaname}{Lemma}
  \providecommand{\propositionname}{Proposition}
\providecommand{\theoremname}{Theorem}
\providecommand{\conjecturename}{Conjecture}
\providecommand{\questionname}{Question}

\title{Hamilton cycles in sparse robustly expanding digraphs}

\author{Allan Lo \thanks{The research leading to these results was partially supported by the  European Research Council under the European Union's Seventh Framework Programme (FP/2007--2013) / ERC Grant Agreement n. 258345 (A. Lo).}\\
\small School of Mathematics\\[-0.8ex]
\small University of Birmingham\\[-0.8ex] 
\small Birmingham, U.K.\\
\small\tt s.a.lo@bham.ac.uk\\
\and
Viresh Patel\thanks{Supported by the Netherlands Organisation for Scientific Research (NWO) through the Gravitation Programme Networks (024.002.003).}\\
\small Korteweg de Vries Instituut voor Wiskunde \\[-0.8ex]
\small Universiteit van Amsterdam\\[-0.8ex]
\small 1090 GE AMSTERDAM, The Netherlands\\
\small\tt  viresh.s.patel@gmail.com}

\begin{document}

\global\long\def\labelenumi{(\roman{enumi})}

\maketitle


\begin{abstract}
The notion of robust expansion has played a central role in the solution of several conjectures involving the packing of Hamilton cycles in graphs and directed graphs. These and other results usually rely on the fact that every robustly expanding (di)graph with suitably large minimum degree contains a Hamilton cycle. Previous proofs of this require Szemer{\'e}di's Regularity Lemma and so this fact can only be applied to dense, sufficiently large robust expanders.
We give a proof that does not use the Regularity Lemma and, indeed, we can apply our result to sparser robustly expanding digraphs.
\end{abstract}

\maketitle

\section{Introduction\label{sec:intro}}

Throughout, we work with simple directed graphs (also called digraphs), i.e.\ directed graphs with no loops and with at most two edges between each pair of vertices (one in each direction).
A Hamilton cycle in a (directed) graph is a (directed) cycle that passes through every vertex. Over the last several decades, there has been intense study in finding sufficient conditions for the existence of Hamilton cycles in graphs and digraphs. The seminal result in the case of graphs is Dirac's Theorem~\cite{Dirac} and in the case of digraphs is Ghouila-Houri's Theorem~\cite{G-H}, each giving tight minimum degree conditions for the existence of Hamilton cycles.

This paper concerns Hamilton cycles in robust expanders. Below we define  a robust expander and give some brief background.

\begin{definition}
For an $n$-vertex digraph $D=(V,E)$, $\nu\in (0,1)$, and $S 
\subseteq V$, the \emph{robust $\nu$-outneighbourhood} of $S$, denoted $\RN^+_{\nu}(S)$, is the set of vertices that have at 
least $\nu n$ inneighbours in $S$. Given $0 < \nu \leq \tau < 1$, 
we say $D$ is a \emph{robust $(\nu,\tau)$-outexpander} if 
\[
|\RN^+_{\nu}(S)| \geq |S| + \nu n
\] 
for every $S \subseteq V$ satisfying $\tau n \leq |S| \leq (1 - \tau)n$. 
The \emph{robust $\nu$-in-neighbourhood}, $\RN^-(S)$, and \emph{robust $(\nu,\tau)$-inexpanders} are defined similarly. We refer to $D$ as a robust $(\nu,\tau)$-expander if it is both a robust $(\nu,\tau)$-in and -outexpander. 
\end{definition}

Usually the parameters $\nu$ and $\tau$ are thought of as small constants as in Theorem~\ref{th:sz} below, but we will also be interested in these parameters as functions of $n$.
Note that robust expansion is a resilient property, i.e. if $D$ is a robust outexpander, then $D$ remains a robust outexpander (with slightly worse parameters) after removing a sparse subgraph. 

Robust expansion has played a central role in the proofs of several conjectures about Hamilton cycles. The starting point of many of these proofs is the following result which says that a robust expander with linear minimum semi-degree contains a Hamilton cycle. The semi-degree $\delta^0(D)$ of a digraph $D$ is given by $\delta^0(D) = \min(\delta^+(D), \delta^-(D))$ where $\delta^+(D)$ and $\delta^-(D)$ are respectively the minimum outdegree and minimum indegree of $D$.

\begin{theorem}[\cite{KOT2}]
\label{th:sz}
Let $n_0$ be a positive integer and $\gamma, \nu, \tau$ be positive constants such
that $1/n_0 \ll \nu \leq \tau \ll \gamma < 1$. 
Let $D$ be a digraph on $n \geq n_0$ vertices with
$\delta^0(D) \geq \gamma n$ which is a robust $(\nu, \tau)$-outexpander. Then $D$ contains a Hamilton
cycle.
\end{theorem}
This result was first proved in~\cite{KOT2} by K\"uhn, Osthus and Treglown. A simpler proof is given in \cite{KOsurvey} and an algorithmic version is given in \cite{CKKO1}.
The proofs of Theorem~\ref{th:sz} presented in \cite{KOT2,KOsurvey,CKKO1} all rely on the Regularity Lemma and so in particular one can only work with sufficiently large and dense digraphs. 

Our main purpose in this paper is to give a proof of the above result that avoids the use of the Regularity Lemma, but uses instead the recent absorption technique developed by R{\"o}dl, Ruci{\'n}ski and Szemer{\'e}di~\cite{RRS} (with special forms appearing in earlier work e.g.\ \cite{Kriv}). We apply our technique to ``sparse'' robust expanders which have not been studied before but which we hope may find applications. In addition we consider cycles of different lengths. The most general form of our result is stated below.

\begin{theorem}
\label{th:maino}
Let $n \in \mathbb{N}$ and $\nu, \tau, \gamma \in(0,1)$ satisfying $4 \sqrt[13]{ \log^2 n /n} < \nu \le \tau \le \gamma /16 < 1/16$.
Let $D$ be an $n$-vertex digraph with $\delta^0(G) \geq \gamma n$ which is a robust $(\nu, \tau)$-outexpander. Then, for any $ \nu n/2 \le \ell \le n$ and any vertex~$v$ of~$D$, $D$ contains a directed cycle of length $\ell$ through~$v$.
\end{theorem}

The result above is algorithmic. We believe some form of it should be true for much sparser graphs than we are able to prove it for.

Theorem~\ref{th:sz} (and its undirected version) have been used as a black box in several papers including \cite{KOT1,KMO1,KT,Kelly,OS,FKS}.
Below we discuss results that require the Regularity Lemma only because they rely (directly or indirectly) on Theorem~\ref{th:sz}. For some such results, we can now replace Theorem~\ref{th:sz} with Theorem~\ref{th:maino} to give proofs that do not require the Regularity Lemma and consequently hold for much smaller values of $n$.

\subsection{Hamiltonicity in oriented graphs}
Robust expansion was first used to prove an approximate analogue of Dirac's Theorem for oriented graphs (an oriented graph is a directed graph in which there is at most one edge between each pair of vertices).
\begin{theorem}[\cite{KelKO}]
\label{thm:luke}
For every $\varepsilon >0$ there exists $n_0 = n_0(\varepsilon)$ such that if $D$ is an oriented graph with $n >n_0$ vertices and $\delta^0(D) > \frac{3}{8}n + \varepsilon n$ then $D$ contains a Hamilton cycle.
\end{theorem}
Here the constant $3/8$ cannot be improved due to examples given in \cite{KelKO}. The result above was proved using the Regularity Lemma and an exact version was proved later in \cite{KeeKO} also using the Regularity Lemma. A consequence of Theorem~\ref{th:maino} is that one can adapt the proof of Theorem~\ref{thm:luke} to avoid the use of the Regularity Lemma.
\begin{corollary}
\label{co:luke}
Let $n \in \mathbb{N}$ and $0<\varepsilon < 1/64$ with $n >  \varepsilon^{-40} $. If $D$ is an $n$-vertex oriented graph with $\delta^0(D) > \frac{3}{8}n + \varepsilon n$ then $D$ contains a Hamilton cycle.
\end{corollary} 
\noindent
In fact, one can use Theorem~\ref{th:maino} to adapt the proof of the exact version in~\cite{KeeKO} to avoid the use of the Regularity Lemma.

\subsection{Hamiltonicity and degree sequences}
In~\cite{KOT2}, K{\"u}hn, Osthus and Treglown give an approximate solution to a conjecture of Nash-Williams~\cite{NW} about sufficient conditions on the degree sequence of a digraph to guarantee the existence of a Hamilton cycle. Their result uses the Regularity Lemma, but Theorem~\ref{th:maino} can be used to adapt their proof to avoid the use of the Regularity Lemma and thus give a better approximation.

For a digraph~$D$, consider its outdegree sequence $d_1^+\le \dots \le d_n^+$ and indegree sequence $d_1^- \le \dots \le d_n^-$.
Note that $d_i^+$ and $d_i^-$ do not necessarily correspond to the degree of the same vertex of~$D$.
\begin{theorem} \label{NWa}
Let $n \in \mathbb{N}$ and $\gamma \in (0,1/2)$ be such that $n \ge 2^{91} \gamma^{-27}$.
Let $D$ be an $n$-vertex digraph such that for all $i < n/2$,
\begin{itemize}
	\item $d_i^+ \ge i + \gamma n$ or $d_{n-i-\gamma n }^- \ge n-i$,
	\item $d_i^- \ge i + \gamma n$ or $d_{n-i-\gamma n }^+ \ge n-i$.
\end{itemize}
Then, for any $ \nu n/2 \le \ell \le n$ and any vertex $v$ of~$D$, $D$ contains a directed cycle of length $\ell$ through~$v$.
\end{theorem}

\subsection{Hamiltonicity in regular graphs}
In~\cite{KLOS1, KLOS2}, K\"uhn, Osthus, Staden and the first author prove the one remaining case of a conjecture of Bollob{\'a}s and Haggvist, making (indirect) use of the Regularity Lemma: they prove that there exits $n_0$ such that every $3$-connected $D$-regular graph on $n \geq n_0$ vertices with $D \ge n/4$ is Hamiltonian. 
Replacing the use of Theorem~\ref{th:sz} by Theorem~\ref{th:maino} in~\cite{KLOS1, KLOS2} gives a proof of the result avoiding the Regularity Lemma.

\subsection{Outline} 

In the next section we collect some notation and in Section~\ref{sec:preliminaries}, we prove some simple facts about robustly expanding digraphs. Section~\ref{sec:absorb} is devoted to describing and constructing an `absorbing structure' $H$ in a robustly expanding digraph $D$. Informally, one can think of $H$ as a set of edges of $D$ which have the property that (almost) any small collection of vertex-disjoint cycles of $D$ can be connected  together into a long cycle using the edges of $H$. In Section~\ref{sec:r-e} we show that the vertices of any robustly expanding digraph can be covered by a small number of cycles. In Section~\ref{sec:Ham} we combine these results to prove Theorem~\ref{th:maino}, and we give some concluding remarks in Section~\ref{sec:conc}.

We mention here that during the course of various proofs, several straightforward calculations, which we feel detract from the main argument, are suppressed and can be found at the end of the paper.

\section{Notation\label{sec:notation}}

The digraphs considered in this paper do not have loops and we allow up to two edges between any pair $x$, $y$ of distinct vertices, at most one in each direction. Given a digraph $D=(V,E)$, we sometimes write $V(D):=V$ for its vertex set and $E(D):=E$ for its edge set and $|D|$ for the number of its vertices. We write $xy$ for an edge directed from $x$ to $y$.

We write $H \subseteq D$ to mean that $H$ is a subdigraph of $D$, i.e.\ $V(H) \subseteq V(D)$ and $E(H) \subseteq E(D)$.
Given $X\subseteq V(D)$, we write $D-X$ for the digraph obtained from $D$ by deleting all vertices in~$X$,
and $D[X]$ for the subdigraph of $D$ induced by~$X$.
Given $F\subseteq E(D)$, we write $D- F$ for the digraph obtained from $D$ by deleting all edges in~$F$. If $H$ is a subdigraph of $D$, we write $D- H$ for $D- E(H)$. For two subdigraphs $H_1$ and $H_2$ of $D$, we write $H_1 \cup H_2$ for the subdigraph with vertex set $V(H_1) \cup V(H_2)$ and edge set $E(H_1) \cup E(H_2)$. 
For a set $U$, $U^2$ means the set of all ordered pairs of $U$, and $U^{[2]}$ means the set of all ordered pairs of $U$ except pairs of the form $(x,x)$.

If $x$ is a vertex of a digraph $D$, then $N^+_D(x)$ denotes the \emph{outneighbourhood} of $x$, i.e.~the
set of all those vertices $y$ for which $xy\in E(D)$. Similarly, $N^-_D(x)$ denotes the \emph{inneighbourhood} of $x$, i.e.~the
set of all those vertices $y$ for which $yx\in E(D)$. 
We write $d^+_D(x):=|N^+_D(x)|$ for the \emph{outdegree} of $x$ and $d^-_D(x):=|N^-_D(x)|$ for its \emph{indegree}. 
We denote the \emph{minimum outdegree} of $D$ by $\delta^+(D):=\min \{d^+_D(x): x\in V(D)\}$
and the \emph{minimum indegree} $\delta^-(D):=\min \{d^-_D(x): x\in V(D)\}$.
The \emph{minimum semi-degree} of $D$ is $\delta^0(D):=\min\{\delta^+(D), \delta^-(D)\}$.

Unless stated otherwise, when we refer to paths and cycles in digraphs, we mean
directed paths and cycles, i.e.~the edges on these paths and cycles are oriented consistently. We write $P = x_1x_2 \cdots x_t$ to indicate that $P$ is a path with edges $x_1x_2, x_2x_3, \ldots, x_{t-1}x_t$, where $x_1, \ldots, x_t$ are distinct vertices. We occasionally denote such a path $P$ by $x_1Px_t$ to indicate that it starts at $x_1$ and ends at $x_t$. 
We write $\mathring{P}$ for the interior of $P$, i.e.\ $\mathring{P} = x_2 \cdots x_{t-1}$.
For two paths $P=a \cdots b$ and $Q = b \cdots c$, we write $aPbQc$ for the concatenation of the paths $P$ and $Q$ and this notation generalises to cycles in the obvious ways. 

Throughout, logarithms are taken base $e$.

\section{Preliminaries\label{sec:preliminaries}}

In this section, we prove some basic properties of robust expanders.
The following proposition follows immediately from the definition of a robust expander.
\begin{proposition}
\label{pr:del}
Suppose $D=(V,E)$ is a robust $( \nu, \tau)$-expander and $S \subseteq V$ with $|S| \leq \varepsilon n$. Then $D - S$ is a $( \nu - \varepsilon, \tau / (1- \varepsilon))$-expander.
\end{proposition}

The following observation of DeBiasio, which can be found in \cite{Tay}, says that robust inexpansion is essentially equivalent to robust outexpansion; thus we can and will restrict ourselves to digraphs that are robust $(\nu,\tau)$-expanders. We reproduce the proof explicitly quantifying the relationships between the various parameters.
\begin{proposition}[DeBiasio] 
\label{pr:dib}
Suppose $D=(V,E)$ is an $n$-vertex robust $(\nu, \tau)$-outexpander with $\delta^0(D) \geq \gamma n$, where $\gamma > 2\tau$, $\tau \gamma > \nu^2/2$ and $\nu < 1/2$. Then $D$ is a robust $(\nu^2/2, 2\tau)$-inexpander. 
\end{proposition}

\begin{proof}
Suppose that $D$ is not a robust $(\nu^2/2, 2\tau)$-inexpander. Then there is a set $S \subseteq V$ with $2\tau n \leq |S| \leq (1- 2\tau)n$ such that $|\RN^-_{\nu^2/2}(S)| < |S| + \nu^2/2$. Let $T = V \setminus RN^-_{\nu^2/2}(S)$. Observe that
\[
|S|\gamma n \leq e(V,S) \leq |RN^-_{\nu^2/2}(S)||S| + \nu^2  n^2/ 2,
\]
so $|RN^-_{\nu^2 / 2}(S)| \geq \gamma n/2$, where we used that $|S|\gamma n /2 \ge \tau \gamma n^2 \ge \nu^2 n^2/2$. Therefore
\[
\tau n < n - (1 - 2\tau + \nu^2/2)n < |T| \le (1 - \gamma/2)n < (1-\tau)n,
\]
where the first and last inequalities follow from our choice of parameters. By the definition of $T$, we have that $e(T,S) < |T|\nu^2 n /2$ and so $|\RN^+_{\nu}(T) \cap S| < |T| \nu /2 < \nu n / 2$. Hence
\begin{align*}
|\RN^+_{\nu}(T)| 
&= |\RN^+_{\nu}(T) \setminus S| + |\RN^+_{\nu}(T) \cap S| \\
&< (n - |S|) + \nu n /2
\le n - (|S| + \nu^2 n / 2 ) + \nu n /4 \\
&< n - |\RN^-_{\nu^2/2}(S)| + \nu n/4 \le |T| + \nu n,
\end{align*}
where we used that $\nu < \frac{1}{2}$ on the second line. Thus $D$ is not a robust $(\nu, \tau)$-outexpander, a contradiction.
\end{proof}

The next two lemmas show that robust expansion allows us to construct short paths between prescribed pairs of vertices.

\begin{lemma}
\label{le:path}
Let $0 < \nu \le \tau \le \gamma/2 <1/2$ and $n \in \mathbb{N}$ satisfying $n \geq 4 \nu^{-2} $.
Suppose that $D$ is an $n$-vertex digraph which is a robust $( \nu, \tau)$-expander 
and $\delta^{0}(D) \geq \gamma n$. Given distinct vertices $u,v \in V(D)$, 
there exists a path $P = x_0 \cdots x_{t+1}$ in $D$ where $x_0=u$, $x_{t+1}=v$ 
and $t \leq \nu^{-1} - 1$. (Note that $P$ consists of at most $\nu^{-1}+1$ vertices.) 
\end{lemma}

\begin{proof}
Let $N_1:=N^+(u)$ and inductively define $N_{i+1} := RN_{\nu}^+(N_i)$.
Note that $|N_1| \geq \gamma n > \tau n$, so for all $i \geq 1$ if $|N_i| < (1- \tau)n$ then $|N_{i+1}| \geq |N_i| + \nu n$.  Observe that, for some $t \leq \nu^{-1}-1$, we have $|N_t| \geq (1 - \tau)n$.
Moreover $N_{t+1} = V(D)$ as $\delta^{0}(D) \geq \gamma n \ge 2 \tau n$.

Set $x_{t+1} = v$.
For $i = t,t-1, \dots,1$, let $x_i$ be a vertex in $( N^-(x_{i+1}) \cap N_i )\setminus \{u,x_{i+1}, \dots, x_{t+1}\}$, which exists as $x_{i+1} \in RN_{\nu}^+(N_i)$ implies that $ | N^-(x_{i+1}) \cap N_i | \geq \nu n \geq \nu^{-1} + 1 \geq t+2$.
By setting $x_0 = u$, we obtain a direct path $P=x_0 \cdots x_{t+1}$ in $D$.
\end{proof}

\begin{lemma}
\label{le:paths}
Let $0 < \nu \le \tau \le \gamma/4 <1/4$ and $n,r \in \mathbb{N}$ satisfying $n \geq ( 6 r + 11) \nu^{-2}$.
Suppose that $D$ is an $n$-vertex digraph which is a robust $(\nu, \tau)$-expander and $\delta^{0}(D) \geq \gamma n$. Given distinct vertices $u_1 \ldots, u_r,v_1, \ldots, v_r \in V(D)$, there exists vertex-disjoint paths $P_1, \ldots, P_r$ in $D$ where $P_i$ is from $u_i$ to $v_i$ and $|P_i| \leq 2\nu^{-1}+1$.
\end{lemma}

\begin{proof}
By induction assume that we have constructed vertex-disjoint paths $P_1, \ldots, P_{k-1}$ in $D$ for some $k < r$, where, for each $i= 1, \ldots, k-1$, $P_i$ is from $u_i$ to $v_i$ and $|P_i| \leq 2\nu^{-1}+1 \leq 3\nu^{-1}$ and $V(P_i) \cap \{u_{i+1}, \ldots, u_r, v_{i+1}, \ldots, v_r \} = \emptyset$. Let $D_{k-1}$ be the digraph obtained from $D$ by deleting all vertices in $P_1, \ldots, P_{k-1}$ and all vertices $u_{k+1}, \ldots, u_r, v_{k+1}, \ldots, v_r$.
Note that $D_{k-1}$ is obtained from $D$ by deleting at most $3 r \nu^{-1} \leq \frac12 \nu  n$ vertices, so by Proposition~\ref{pr:del}, $D_{k-1}$ is a robust $( \frac{1}{2}\nu, \frac{8}{7} \tau)$-expander with $\delta^{0}(D_{k-1}) \geq \frac{7}{8}\gamma n \geq \frac{7}{8}\gamma |D_{k-1}|$.
Note that $|D_{k-1}| \geq n - 3 r \nu^{-1} \geq n -  r \nu^{-2} \geq 16  \nu^{-2}$.
Apply Lemma~\ref{le:path} to $D_{k-1}$ giving a path $P_k$ in $D_{k-1}$ of length at most $2\nu^{-1} +1$ from $u_k$ to $v_k$. Thus $P_k$ (as a path in $D$) is vertex-disjoint from $P_1, \ldots, P_{k-1}$ and $\{u_{k+1}, \ldots, u_r, v_{k+1}, \ldots, v_r\}$ as required. Thus by induction we can find the paths $P_1, \ldots, P_{r}$.
\end{proof}

We give a simple inequality that we shall use several times.

\begin{proposition}
\label{pr:polyexp}
Fix $k, a > 0$. Then $e^x>ax^k$ for all $x \geq \max(3k (\log k +1) + 3\log a,0)$. Similarly for $c,d>0$ if $x> 3c(\log c +1) + 3d$, then we have $x > c \log x + d$.
\end{proposition}
\begin{proof}
We start by showing that for all $a>0$ and $x \geq \max(3\log a + 3, 0)$, we have $e^x \geq ax$. This is clearly true if $0 < a \leq 1$. If $a>1$, set $f(x) = e^x - ax$ and set $x_0 := 3\log a + 3>0$. We have $f(x_0) = e^3a^3 - 3a\log a - 3a>0$ and $f'(x) = e^x - a >0$ for all $x \geq x_0$. Hence $f(x)>0$ for all $x \geq x_0$ and so $e^x > ax$ for all $x \geq \max(3\log a + 3,0)$.

Finally, making the transformation $X = kx$ and $A=a^k/k^k$, and assuming $A,k > 0$, the inequality above becomes $e^X \geq AX^k$ for all $X \geq \max(3k\log k + 3\log A + 3k, 0)$.

For the other inequality, note that $x > c\log x + d$ if and only if $e^x > e^dx^c$, which holds if $x> \max(3c\log c + 3 d + 3c,0)$.
\end{proof}


\section{The absorbing structure\label{sec:absorb}}

In this section, we describe what we mean by an \emph{absorbing structure} and show how to find one in a robustly expanding digraph with large minimum in- and outdegree.
We begin by informally describing the properties we desire our absorbing structure to have. Given a digraph $D$ we shall seek a subdigraph $S \subseteq D$ with the properties that 
\begin{itemize}
\item $|V(S)|$ is small;
\item $S$ contains a Hamilton cycle (on $V(S)$);
\item In $D$, given a small number of any vertex-disjoint paths $P_1, \ldots, P_d$ that are also vertex-disjoint from $S$, we can use $S$ to absorb $P_1, \ldots, P_d$ into $C$ i.e.\ we can find a Hamilton cycle $C'$ on $V(S) \cup (\bigcup_{i=1}^d V(P_i) )$.
\end{itemize}
The sequence of definitions that follow will lead to a precise description of our absorbing structure. We start by defining an alternating path.

\begin{definition}
Let $D$ be a digraph, and let $x_1, \ldots, x_t$ be distinct
 vertices of $D$ with $t$ even. An \emph{alternating path} 
$P = [x_1x_2 \cdots x_t]$ is a subgraph of $D$ with vertex
 set $\{x_1, \ldots, x_t\}$ (where $x_1, \ldots, x_t$ are distinct vertices)  and edge set 
\[
\{x_ix_{i+1} \mid i=1,3,5, \ldots, t-1 \} \cup \{x_{j+1}x_j \mid j=2,4,6, \ldots, t-2 \}.
\]
We say $P$ is an alternating path from $x_1$ to $x_t$.
\end{definition}
An alternating path is thus a path where the directions of the edges alternate. It will be important for us that the number of vertices in an alternating path is even so that the first vertex has outdegree $1$ and the last vertex has indegree $1$. 

As with paths, robust expansion allows us to construct alternating paths between prescribed vertices.
\begin{lemma}
\label{le:altpath}
Let $0 < \nu \le \tau \le \gamma/2 <1/2$ and $n \in \mathbb{N}$ satisfying $n \geq 4 \nu^{-2}$.
Suppose that $D$ is an $n$-vertex digraph which is a robust $(\nu, \tau)$-expander and $\delta^{0}(D) \geq \gamma n$. Given distinct vertices $u,v \in V(D)$, there exists an alternating path $P = [x_0 \cdots x_tx_t^* \cdots x_0^*]$ in $D$ where $x_0=u$, $x_0^*=v$ and $t \leq (\nu^{-1}+4)/2$ is even. (Thus $P$ consists of at most $\nu^{-1}+6$ vertices.) 
\end{lemma}
\begin{proof}
Let $N_1:=N^+(u)$ and inductively define
\begin{align*}
N_{i+1} := 
\begin{cases}
\RN_{\nu}^+(N_i) &\text{ if } i \text{ even;} \\
\RN_{\nu}^-(N_i) &\text{ if } i \text{ odd.}
\end{cases}
\end{align*}
Note that $|N_1| \geq \gamma n > \tau n$, so for all $i \geq 1$ if $|N_i| < (1- \tau)n$ then $|N_{i+1}| \geq |N_i| + \nu n$. Set $r:= \lceil \nu^{-1} \rceil$ and observe that $|N_r| \geq (1 - \tau)n$. 
Note that $N_{r''} = V(D)$ for all $r'' > r$ as $\delta^{0}(D) \geq \gamma n\ge 2 \nu n$.
Choose $r'$ to be the smallest integer that is greater than $r$ and divisible by $4$; thus $r' \leq \nu^{-1} + 4$ and $N_{r'} = V(D)$.

Set $y_{r'+1} = v$.
For $i = r',r'-1, \dots,1$, let $y_i$ be a vertex such that 
\[
y_{i} \in 
\begin{cases}
N^-(y_{i+1}) \cap N_i \setminus \{u,y_{i+1}, \dots, y_{r'+1}\} &\text{if } i \text{ even}; \\
N^+(y_{i+1}) \cap N_i \setminus \{u,y_{i+1}, \dots, y_{r'+1}\} &\text{if } i \text{ odd}.
\end{cases}
\]
To see that such $y_i$ exists, observe that since $y_{i+1} \in N_{i+1} = \RN_{\nu}^+(N_i)$, $|N^-(y_{i+1}) \cap N_i| \ge  \nu n \geq 4\nu^{-1}  \geq \nu^{-1}+ 6 \geq r'+2 $ if $i$ is even (and a similar inequality holds if $i$ is odd). 

Thus we obtain distinct vertices $y_1, \ldots, y_{r'}$ such that $y_{i+1}y_{i},y_jy_{j+1} \in E(D)$ for $i = 1,3,5 \ldots, r'-1$ and $j=2,4,6, \ldots, r'-2$. Then relabelling $y_1, \ldots, y_{r'}$ to $x_1, \ldots, x_tx^*_t, \ldots, x^*_1$ respectively and $x_0:=u, x_0^*:=v$ gives the desired alternating path. Since $r'$ is divisible by $4$, we have that $t$ is even as required.
\end{proof}

Next we define ladders, which will be the key structures that allow us to absorb paths.
\begin{definition}
Let $D$ be a digraph and let $u,v \in V(D)$ be distinct vertices. A \emph{ladder} $L$ from $u$ to $v$ is a subdigraph of $D$ given by 
\[
L = Q \cup Q_1 \cup Q_3 \cup Q_5 \cup \cdots \cup Q_{t-1},
\]
where
\begin{itemize}
\item[(i)] $Q = [x_0 x_1 \cdots x_t x_t^* \cdots x_1^* x_0^*]$ is an alternating path (with $t$ even) and $x_0=u$ and $x_0^*=v$;
\item[(ii)] $Q_i$ is a directed path from $x_i$ to $x_i^*$ for each $i=1,3,\ldots, t-1$; and
\item[(iii)] $Q_1, Q_3, \ldots, Q_{t-1}$ are vertex-disjoint paths and are each internally vertex-disjoint from $Q$.
\end{itemize}
We call $Q$ the alternating path of $L$.
\begin{itemize}
\item For $i=0,2,4, \ldots, t-2$, we define $R_i \subseteq L$ to be the path $R_i := x_ix_{i+1}Q_{i+1}x^*_{i+1}x^*_i$ and $R_t:=x_tx_t^*$. We call these the \emph{rung paths} of $L$.
\item For $i=2,4, \ldots, t$, define $R_i' \subseteq L$ to be the path $R'_i := x_ix_{i-1}Q_{i-1}x^*_{i-1}x^*_i$. We call these the \emph{alternative rung paths} of $L$.
\end{itemize}
We say the ladder $L$ is \emph{embedded} in the cycle $C$ if $R_i \subseteq C$ for all $i=0,2,4, \ldots, t$. 
\end{definition}

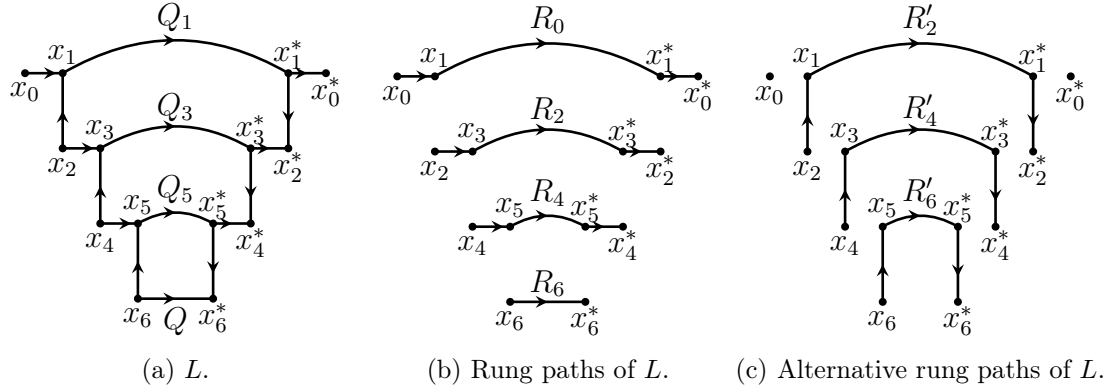
\begin{figure}[t]
\centering
\subfloat[$L$.]{
\begin{tikzpicture}[scale =0.5, line width=1pt]
	\tikzset{middlearrow/.style={
        decoration={markings, mark= at position 0.5 with {\arrow{#1}} ,},
        postaction={decorate}}
	}
     		\foreach \x in {0,1,2}
			{
			\filldraw[fill=black] (\x+1,2*\x) circle (2pt);
			\filldraw[fill=black] (\x+1, 2*\x+2) circle (2pt);
			\filldraw[fill=black] (-\x-1, 2*\x) circle (2pt);
			\filldraw[fill=black] (-\x-1, 2*\x+2) circle (2pt);
			\draw[middlearrow={stealth reversed}] (\x+1 , 2*\x)	-- (\x+1,2*\x+2); 
			\draw[middlearrow={stealth}] (\x+1,2*\x+2)-- (\x+2,2*\x+2);
			\draw[middlearrow={stealth}] (-\x-1, 2*\x)	-- (-\x-1, 2*\x+2);
			\draw[middlearrow={stealth reversed}] (-\x-1,2*\x+2)-- (-\x-2,2*\x+2);
			}
		\filldraw[fill=black] (-4,6) circle (2pt);
		\filldraw[fill=black] (4,6) circle (2pt);
		\draw[middlearrow={stealth}] (-1,0)	-- (1,0);
		
		\draw[middlearrow={stealth}] (-1,2) to[out = 30,in = 150] (1,2);
		\draw[middlearrow={stealth}] (-2,4) to[out = 30,in = 150] (2,4);
		\draw[middlearrow={stealth}] (-3,6) to[out = 30,in = 150] (3,6);

		\node  at (-4,5.5) {$x_0$};
		\node  at (-3,6.5) {$x_1$};
		\node  at (-3,3.5) {$x_2$};
		\node  at (-2,4.5) {$x_3$};
		\node  at (-2,1.5) {$x_4$};
		\node  at (-1,2.5) {$x_5$};
		\node  at (-1,-0.5) {$x_6$};

		\node  at (4,5.5) {$x_0^*$};
		\node  at (3,6.5) {$x_1^*$};
		\node  at (3,3.5) {$x_2^*$};
		\node  at (2,4.5) {$x_3^*$};
		\node  at (2,1.5) {$x_4^*$};
		\node  at (1,2.5) {$x_5^*$};
		\node  at (1,-0.5) {$x_6^*$};
		
		\node  at (0,2.9) {$Q_5$};
		\node  at (0,5.1) {$Q_3$};
		\node  at (0,7.5) {$Q_1$};
		\node  at (0,-0.6) {$Q$};
\end{tikzpicture}
}
\subfloat[Rung paths of $L$.]{
\begin{tikzpicture}[scale =0.5, line width=1pt]
	\tikzset{middlearrow/.style={
        decoration={markings, mark= at position 0.5 with {\arrow{#1}} ,},
        postaction={decorate}}
	}
     		\foreach \x in {0,1,2}
			{
			\filldraw[fill=black] (\x+1,2*\x) circle (2pt);
			\filldraw[fill=black] (\x+1, 2*\x+2) circle (2pt);
			\filldraw[fill=black] (-\x-1, 2*\x) circle (2pt);
			\filldraw[fill=black] (-\x-1, 2*\x+2) circle (2pt);
			\draw[middlearrow={stealth}] (\x+1,2*\x+2)-- (\x+2,2*\x+2);
			\draw[middlearrow={stealth reversed}] (-\x-1,2*\x+2)-- (-\x-2,2*\x+2);
			}
		\filldraw[fill=black] (-4,6) circle (2pt);
		\filldraw[fill=black] (4,6) circle (2pt);
		\draw[middlearrow={stealth}] (-1,0)	-- (1,0);
		
		\draw[middlearrow={stealth}] (-1,2) to[out = 30,in = 150] (1,2);
		\draw[middlearrow={stealth}] (-2,4) to[out = 30,in = 150] (2,4);
		\draw[middlearrow={stealth}] (-3,6) to[out = 30,in = 150] (3,6);

		\node  at (-4,5.5) {$x_0$};
		\node  at (-3,6.5) {$x_1$};
		\node  at (-3,3.5) {$x_2$};
		\node  at (-2,4.5) {$x_3$};
		\node  at (-2,1.5) {$x_4$};
		\node  at (-1,2.5) {$x_5$};
		\node  at (-1,-0.5) {$x_6$};

		\node  at (4,5.5) {$x_0^*$};
		\node  at (3,6.5) {$x_1^*$};
		\node  at (3,3.5) {$x_2^*$};
		\node  at (2,4.5) {$x_3^*$};
		\node  at (2,1.5) {$x_4^*$};
		\node  at (1,2.5) {$x_5^*$};
		\node  at (1,-0.5) {$x_6^*$};
		
		\node  at (0,2.9) {$R_4$};
		\node  at (0,5.1) {$R_2$};
		\node  at (0,7.5) {$R_0$};
		\node  at (0,0.5) {$R_6$};
\end{tikzpicture}
}
\subfloat[Alternative rung paths of $L$.]{
\begin{tikzpicture}[scale =0.5, line width=1pt]
	\tikzset{middlearrow/.style={
        decoration={markings, mark= at position 0.5 with {\arrow{#1}} ,},
        postaction={decorate}}
	}
     		\foreach \x in {0,1,2}
			{
			\filldraw[fill=black] (\x+1,2*\x) circle (2pt);
			\filldraw[fill=black] (\x+1, 2*\x+2) circle (2pt);
			\filldraw[fill=black] (-\x-1, 2*\x) circle (2pt);
			\filldraw[fill=black] (-\x-1, 2*\x+2) circle (2pt);
			\draw[middlearrow={stealth reversed}] (\x+1 , 2*\x)	-- (\x+1,2*\x+2); 
			\draw[middlearrow={stealth}] (-\x-1, 2*\x)	-- (-\x-1, 2*\x+2);
			}
		\filldraw[fill=black] (-4,6) circle (2pt);
		\filldraw[fill=black] (4,6) circle (2pt);
		
		\draw[middlearrow={stealth}] (-1,2) to[out = 30,in = 150] (1,2);
		\draw[middlearrow={stealth}] (-2,4) to[out = 30,in = 150] (2,4);
		\draw[middlearrow={stealth}] (-3,6) to[out = 30,in = 150] (3,6);

		\node  at (-4,5.5) {$x_0$};
		\node  at (-3,6.5) {$x_1$};
		\node  at (-3,3.5) {$x_2$};
		\node  at (-2,4.5) {$x_3$};
		\node  at (-2,1.5) {$x_4$};
		\node  at (-1,2.5) {$x_5$};
		\node  at (-1,-0.5) {$x_6$};

		\node  at (4,5.5) {$x_0^*$};
		\node  at (3,6.5) {$x_1^*$};
		\node  at (3,3.5) {$x_2^*$};
		\node  at (2,4.5) {$x_3^*$};
		\node  at (2,1.5) {$x_4^*$};
		\node  at (1,2.5) {$x_5^*$};
		\node  at (1,-0.5) {$x_6^*$};
		
		\node  at (0,2.9) {$R'_6$};
		\node  at (0,5.1) {$R'_4$};
		\node  at (0,7.5) {$R'_2$};
\end{tikzpicture}
}
\caption{A Ladder $L = Q  \cup Q_1 \cup Q_3 \cup Q_5$}
\label{fig:ladder}
\end{figure}

It is relatively easy to construct ladders in robust expanders.
First we show how a ladder embedded in a cycle can be used to absorb a path into the cycle.

\begin{lemma}
\label{le:1-absorb}
Let $D$ be a digraph and let $u,v \in V(D)$ be distinct vertices. Let $L \subseteq D$ be a ladder from $u$ to $v$ embedded in a cycle $C \subseteq D$. For any path $P \subseteq D$ from $u$ to $v$ that is internally vertex-disjoint from $C$ there exists a cycle $C' \subseteq D$ such that
\begin{itemize}
\item[(i)] $P \subseteq C'$, 
\item[(ii)] $V(L) \subseteq V(C')$, 
\item[(iii)] for any path $P' \subseteq C$ with $V(P') \cap V(L) = \emptyset$, we have $P' \subseteq C'$, and
\item[(iv)] if $x \in V(D)$ satisfies $x \not\in V(C) \cup V(P)$, then $x \not\in V(C')$.
\end{itemize} 
In particular, (i), (ii), (iii), and (iv) immediately imply
\begin{itemize}
\item [(v)] $V(C) \cup V(P) = V(C')$.
\end{itemize}
\end{lemma}
\begin{proof}
Let $Q = [x_0 \cdots x_tx_t^* \cdots x_0^*]$ be the alternating path of $L$, and let $Q_i$ be the corresponding paths of $L$ from $x_i$ to $x_i^*$ for $i = 1,3, \ldots, t-1$.
 Let $R_0, R_2, \ldots, R_t$ be the rung paths of $L$ and  $R_2', R_4', \ldots, R_t'$ the alternative rung paths. Set $R_0':=P$. We simply replace $R_i$ with $R_i'$ in $C$ one at a time to obtain the desired cycle $C'$. We spell out the details of the induction below.

We define cycles $C_{-2}, C_0, C_2, C_4, \ldots, C_t$ as follows. Set $C_{-2}:=C$. By induction, we assume that $C_{i-2}$ is a cycle with $R_0', \ldots, R_{i-2}', R_{i}, \ldots, R_t \subseteq C_{i-2}$ (implicitly noting these paths are vertex-disjoint) and that $\mathring{R}'_i = Q_{i-1}$ is vertex-disjoint from $C_{i-2}$. We obtain $C_{i}$ by deleting $R_{i}$ from $C_{i-2}$ and replacing it with $R'_{i}$. Since $R_{i}$ and $R'_{i}$ are internally vertex-disjoint and both are paths from $x_i$ to $x_i^*$, then $C_i$ is a cycle. Clearly we have $R_0', \ldots, R_{i}', R_{i+2}, \ldots, R_t \subseteq C_i$. Since $\mathring{R}_i = \mathring{R}'_{i+2} = Q_{i+1}$ is vertex-disjoint from $C_i$ (since we deleted $\mathring{R}_i$), then $\mathring{R}'_{i+2}$ is vertex-disjoint from $C_i$. 

Thus by induction, we have that $C':=C_t$ is a cycle with $R_0', \ldots, R_t' \subseteq C'$. Therefore $P = R_0' \subseteq C'$ proving (i). Furthermore, since
\[
V(\bigcup_{ \substack{i=0 \\ i \text{ even}}}^t R_i ) =  
V(\bigcup_{ \substack{i=2 \\ i \text{ even}}}^t R'_i ) =
V(L)
\]  
then $V(L) \subseteq V(C')$ proving (ii). In the above induction, we note that if $P' \subseteq C_{i-2}$ is a path vertex-disjoint from $L$, then $P' \subseteq C_i$, so by induction if $P' \subseteq C=C_0$ is a path vertex-disjoint from $L$, then $P' \subseteq C_t=C'$ proving (iii). Finally, we note that, in the above induction, for any vertex $x \in V(D) \setminus (V(P) \cup V(L))$, if $x \not \in C_{i-2}$ then $x \not\in C_i$, proving (iv) and completing the proof.
\end{proof}
From the previous lemma, we now see that embedding several carefully chosen ladders into a cycle can give us the absorbing structure we desire. The next definition makes precise what we mean by `carefully' in the previous sentence.

\begin{definition}
Given a digraph $D$ and distinct vertices $x,y,u,v \in V(D)$, we say that the ordered pair $(u,v) \in V(D)^{[2]}$ \emph{covers} $(x,y) \in V(D)^{2}$ if $ux,yv \in E(D)$. Given $K \subseteq V(D)^{[2]}$ and $U \subseteq V(D)$, we say that $K$ \emph{$d$-covers} $U$ if for every $(x,y) \in U^{2}$ there exist $d$ distinct elements of $K$ each of which covers $(x,y)$. We say $K$ is vertex-disjoint if no two elements of $K$ share a vertex.
\end{definition}

Our motivation for this definition is the following. Suppose $L$ is a ladder from $u$ to $v$ embedded in a cycle $C$ and $P$ is a path from $x$ to $y$ that is vertex-disjoint from $C$, and suppose further that $(u,v)$ covers $(x,y)$. Then we can extend $P$ to the path $uxPyv$ and use the previous lemma to absorb $P$ into $C$. For a digraph $D$, if we can find a small set $K \subseteq V(D)^{[2]}$ which $d$-covers $V(D)$, then we might hope to construct vertex-disjoint ladders from $u$ to $v$ for each $(u,v) \in K$ and embed all those ladders into a cycle $C$. This structure would then have the property that any $d$ vertex-disjoint paths of $D$ (that are also vertex-disjoint from $C$) could be absorbed into $C$. This will be our absorbing structure.

\begin{definition}
Given a digraph $D$ and $d \in \mathbb{N}$, a \emph{$d$-absorber} $S$ of $D$ is a triple $S =(K, \mathcal{L}, C)$, where
\begin{itemize}
\item $K \subseteq V(D)^{[2]}$ is a set of vertex-disjoint pairs which $d$-covers $V(D)$, 
\item $\mathcal{L}$ is a set of vertex-disjoint ladders such that for each $(u,v) \in K$, we have a ladder $L \in \mathcal{L}$ from $u$ to $v$, 
\item $C \subseteq D$ is a cycle such that each $L \in \mathcal{L}$ is embedded in $C$. 
\end{itemize}
We sometimes abuse notation by also writing $S$ for the subgraph $(\cup_{L \in \mathcal{L}} L) \cup C$ of $D$. Note that $V(C) = V(S)$.
\end{definition}

It follows from Lemma~\ref{le:1-absorb} that a $d$-absorber can absorb $d$ vertex-disjoint paths into its cycle.
\begin{corollary}
\label{co:d-absorb}
Let $D$ be a digraph and let $S \subseteq D$ be a $d$-absorber. Suppose $P_1, \ldots, P_r$ are vertex-disjoint paths in $D$ that are also vertex-disjoint from $V(S)$ and $r \leq d$. Then there exists a cycle $C^*$ in $D$ such that $V(C^*)=V(S) \cup V(P_1) \cup \cdots \cup V(P_r)$. 
\end{corollary}
\begin{proof}
Let $x_i$ and $y_i$ be such that $P_i$ is a path from $x_i$ to $y_i$ for $i= 1, \ldots, r$ and let $S = (K, \mathcal{L}, C)$. Since $S$ is a $d$-absorber, for each $i= 1, \ldots, r$, there exists $(u_i,v_i) \in K$ and $L_i \in \mathcal{L}$ such that $(u_i,v_i)$ covers $(x_i,y_i)$ and $L_i$ is a ladder from $u_i$ to $v_i$, and where $u_1, \ldots, u_r, v_1, \ldots, v_r$ are distinct vertices. For each $i$, observe that $Q_i := u_ix_iP_iy_iv_i$ is a path in $D$ and that $Q_1, \ldots, Q_r$ are vertex-disjoint.

Set $C_0:=C$ and assume by induction that there is a cycle $C_{i-1} \subseteq D$ with the property that $V(C_{i-1}) = V(C_0) \cup V(Q_1) \cup \cdots \cup V(Q_{i-1})$ and where $L_i, \ldots, L_r$ are embedded in $C_{i-1}$. Since $L_i$ is a ladder from $u_i$ to $v_i$ embedded in $C_{i-1}$ and $Q_i$ is a path from $u_i$ to $v_i$ internally vertex-disjoint from $C_{i-1}$, Lemma~\ref{le:1-absorb} implies that there exists a cycle $C_i$ such that $V(C_i) = V(C_{i-1}) \cup V(Q_i) = V(C) \cup V(Q_1) \cup \cdots \cup V(Q_i)$. Furthermore, by Lemma~\ref{le:1-absorb}, since $L_{i+1}, \ldots, L_r$ are vertex-disjoint from $L_i$, and are embedded in $C_{i-1}$, so they are embedded in $C_i$. 

This completes the induction step and so we obtain a cycle $C^* := C_r$ of $D$ where $V(C_r) = V(C) \cup V(Q_1) \cup \cdots \cup V(Q_r) = V(S) \cup V(P_1) \cup \cdots \cup V(P_r)$.   
\end{proof}

The sequence of lemmas that follow show how to build a $d$-absorber in a robust expander.
The first lemma shows how to find a $d$-cover in a digraph. 

\begin{lemma}
\label{le:cover}
Let $\gamma \in (0,1)$ and $n,d \in \mathbb{N}$ with $d \geq 8$ and 
$$n>10^5 d^2\gamma^{-4}\log^2(100 d\gamma^{-2}).$$
If $D$ is an $n$-vertex digraph with $\delta^{0}(D) \geq \gamma n$ and $U \subseteq V(D)$, then there exists a vertex-disjoint $K \subseteq V(D)^{[2]}$ with 
$|K| = \lceil 24 \gamma^{-2} ( d \log( 24 d\gamma^{-2} + 2\log n )  \rceil $ which $d$-covers $U$. 
\end{lemma}
\begin{proof}
Set 
$m:= \lceil 24 \gamma^{-2} ( d \log( 24 d\gamma^{-2}) + 2\log n ) \rceil $ and construct $K^*$ randomly by taking a set of $m$ elements, each picked independently and uniformly at random, from $V(D)^{[2]}$; thus $K^*$ may not be vertex-disjoint. We have that
\begin{align*}
\mathbb{P}(K^* \text{ is vertex-disjoint}) &= \prod_{i=0}^{m-1} \left( \binom{n-2i}{2} / \left( \binom{n}{2} - i \right) \right) \\
& \geq 
\prod_{i=0}^{m-1} \left( \binom{n-2i}{2} / \left( n^2/2 \right) \right) 
 = \prod_{i=0}^{2m-1}\left(1 - \frac{i}{n}\right) \\
& \geq 1 - \sum_{i=1}^{2m-1}\frac{i}{n}
 \geq 1 - \frac{2 m^2}{n} > \frac12,
\end{align*}
our choice of $m$ and $n$ and applying Proposition~\ref{pr:polyexp}
	\endnote{
	Let $c'= 100 \gamma^{-2}$ and $d' = 50 d \gamma^{-2} \log (24 d \gamma^{-2})$.
	Since $d \ge 8$, we have 
	\begin{align*}
		\sqrt{n}  \ge 300 d \gamma^{-2}\log(100 d\gamma^{-2}) 
		\ge 300 \gamma^{-2} (\log(100 \gamma^{-2}) +1) + 150 d \gamma^{-2} \log (24 d \gamma^{-2}) \ge 3c'(\log c' +1) +3d'.
	\end{align*}
	By Proposition~\ref{pr:polyexp}, we have 
	\begin{align*}
		\sqrt{n} \ge c' \log \sqrt n + d' 
		\ge 50 \gamma^{-2} \log n + 50 d \gamma^{-2} \log (24 d \gamma^{-2}) > 2m. 
	\end{align*}}.

For $(x,y) \in U^{2}$, let ${\rm cov}(x,y)$ be the set of elements in $V(D)^{[2]}$ that cover $(x,y)$. For a uniformly random element $(u,v)$ of $V(D)^{[2]}$, set
\begin{align*}
p:=\mathbb{P}((u,v) \in {\rm cov}(x,y)) \geq \frac{\gamma n(\gamma n - 1)}{n(n-1)} \geq \frac{\gamma^2}{2}, 
\end{align*} 
where the last inequality follows by our choice of $n$. Let $E_{x,y}$ be the number of distinct elements of $K^*$ that cover $(x,y)$ so that $E_{x,y} \sim {\rm bin}(m,p)$. In particular, 
\begin{align*}
\mathbb{P}(E_{x,y} < d) 
& = \sum_{i=0}^{d-1} \binom{m}{i}p^i(1-p)^{m-i} 
\leq (1-p)^{m-d} \sum_{i=0}^{d-1} \frac{m^i}{i!} \\
&\leq  em^d (1-p)^{m-d}
\leq em^d(1 - \gamma^2/2)^{m-d}\\
& \leq m^d \exp( - \gamma^2(m-d)/2) 
\leq m^d \exp(-\gamma^2m/8)
\leq n^{-2}/2.
\end{align*}
by our choices of $m$ and applying Proposition~\ref{pr:polyexp}%
	\endnote{Note $n^2 m^d\exp(-\gamma^2m/8) \le 1/2$ holds if $e^m \ge (2 m^d n^2)^{8 \gamma^{-2}} $, which holds by Proposition~\ref{pr:polyexp} if $m \ge 24 \gamma^{-2} ( d \log( 8 d\gamma^{-2}) + d + \log (2n^2))$, 
	which holds if $ m  \ge 24 \gamma^{-2} ( d\log( 24 d\gamma^{-2}) + 2\log n)$ as $d\log 3 \ge d + \log 2$.}.
Let $X$ be the number of elements of $U^{2}$ not $d$-covered by $K^*$. Then
\begin{align*}
\mathbb{P}(X \geq 0) & \leq \mathbb{E}(X)  
= \sum_{(x,y) \in U^{2}}\mathbb{P}(E_{x,y} < d)  \leq 1/2.
\end{align*}
Therefore $\mathbb{P}(X=0 \text{ and } K^* \text{ is vertex-disjoint}) >0$. 
\end{proof}

Next we show how to build a ladder in a robust expander.
\begin{lemma}
\label{le:laddercon}
Let $0 < \nu \le \tau < \gamma/8 <1/8$ and $n \in \mathbb{N}$ satisfying $n \geq 57 \nu^{-3}$.
Let $D$ be a robust $(\nu, \tau)$-expander on $n$ vertices with $\delta^{0}(D) \geq \gamma n$ and let
$u,v$ be distinct vertices of $D$. Then there exists a ladder $L$ from $u$ to $v$ with $|L| \leq 3 \nu^{-2}$ and where the alternating path of $L$ has at most $2 \nu^{-1}$ vertices.
\end{lemma}
\begin{proof}
By Lemma~\ref{le:altpath}, we can find an alternating path $Q = [x_0 \cdots x_tx_t^* \cdots x_0^*]$, where $x_0=u$, $x_0^* = v$, and $t \leq (\nu^{-1} + 4)/2$ is even (so this alternating path has at most $\nu^{-1}+6 \le 2 \nu^{-1}$ vertices). Next, as in the definition of ladders, we construct vertex-disjoint paths $Q_1, Q_3, \ldots, Q_{t-1}$, where $Q_i$ is from $x_i$ to $x_i^*$ and is vertex-disjoint from $P$ (except at its end points). We do this using Lemma~\ref{le:paths}.

Let $D'$ be the digraph obtained from $D$ by deleting $x_i$ and $x_i^*$ for each even value $i = 0, \ldots, t$; thus we delete $t+2 \leq (\nu^{-1} + 8)/2 \le \nu^{-1}$ vertices and by our choice of large $n$, Proposition~\ref{pr:del} implies\endnote{The choice of $n$ implies $\nu^{-1} \leq \frac{1}{2} \nu n \le \frac{1}{16} \gamma n \le \frac{1}{16} n$.}
that $D'$ is a robust $(\frac{1}{2}\nu, \frac{16}{15}\tau)$-expander with $\delta^{0}(D') \geq \frac{15}{16}\gamma n$. By our choice of parameters and sufficiently large $n$, we can apply Lemma~\ref{le:paths}\endnote{Note that 
\begin{align*}
n - (t+2) \ge n - \nu^{-1} \ge 56 \nu^{-3} \ge 4( 12\nu^{-1}+11) \nu^{-2}
\end{align*}
Thus conditions of Lemma~\ref{le:paths} hold (with $r = t/2 \leq \nu^{-1}/2$ and $ \nu, \tau, \gamma, n$ replaced by $\frac{1}{2}\nu, \frac{16}{15} \tau,  \frac{15}{16}\gamma, n - (t+2)$).}
with $r = t/2$ to obtain vertex-disjoint  paths $Q_1, Q_3, \ldots, Q_{t-1}$ in $D'$ with each $Q_i$ from $x_i$ to $x_i^*$ and of length at most $4\nu^{-1} + 3$. As paths in $D$, these paths are also vertex-disjoint from $Q$ except at their end-points.

Thus the union of the alternating path $Q$ with the paths $Q_1, Q_3, \ldots, Q_{t-1}$ gives a ladder $L$ from $u$ to $v$. We have $|Q| \leq \nu^{-1}+6 \leq  \nu^{-2}/2$, $|Q_i| \leq 4\nu^{-1} + 3 \le 9 \nu^{-1}/2$ for each odd $i < t \leq \nu^{-1}$. Thus $|L| \leq 9 \nu^{-2}/4 +|Q| \leq 3 \nu^{-2} $.
\end{proof}

Next we show that we can build several ladders (between prescribed vertices) in a robustly expanding digraph (for a suitable choice of parameters).
\begin{lemma}
\label{le:ladderscon}
Let $0 < \nu \le \tau \le \gamma/16 < 1/16$ and $n,k \in \mathbb{N}$ satisfying 
$n \ge 460 k  \nu^{-3}$.
Let $D$ be a robust $(\nu, \tau)$-expander on $n$ vertices with $\delta^{0}(D) \geq \gamma n$ and let $u_1, \ldots, u_k,v_1, \ldots, v_k$ be distinct vertices of $D$. Then we can construct vertex-disjoint ladders $L_1, \ldots, L_k$ from $u_i$ to $v_i$ such that $|L_i| \leq 12 \nu^{-2}$ and $|P_i| \leq 8 \nu^{-1}$, where $P_i$ is the alternating path of $L_i$.
\end{lemma}
\begin{proof}
By induction, suppose we have constructed vertex-disjoint ladders $L_1, \ldots, L_{i-1}$ for some $i \leq k$ where $L_j$ is from $u_j$ to $v_j$ and $|L_j| \leq 12 \nu^{-2}$ for all $j < i$,  where the alternating path $P_j$ of $L_j$ satisfies $|P_j| \leq 8\nu^{-1}$ for all $j < i$, and where $L_1, \ldots, L_{i-1}$ are disjoint from $S_i := \{u_i, \ldots, u_k, v_i, \ldots, v_k \}$. Let $D_i$ be obtained from $D$ by deleting all the vertices of $L_1, \ldots, L_{i-1}$ and $S_{i}\setminus \{u_i,v_i\}$ (so $u_i,v_i \in V(D_i)$); thus the number of vertices deleted is at most
\[
12 \nu^{-2} (i-1) + 2(k-i) \leq 12 k \nu^{-2} \leq \nu n/2 ,
\]
where the last inequality follows from our choice of $n$.
 The inequality above together with Proposition~\ref{pr:del} implies that $D_i$ is a robust $( \frac{1}{2}\nu, \frac{32}{31}\tau)$-expander with $\delta(D_i) \geq \frac{31}{32}\gamma n$. By our choice of parameters and $n$, we can apply\endnote{We check the conditions of Lemma~\ref{le:laddercon} with $ \nu, \tau, \gamma, n$ replaced by $ \frac{1}{2}\nu, \frac{32}{31}\tau, \frac{31}{32}\gamma,|D_i|$.
Note that $	|D_i| \ge n -  12 k \nu^{-2} \ge 459 \nu^{-3} \ge 57 (\nu/2)^{-3} $.} Lemma~\ref{le:laddercon}
to obtain a ladder $L_i$ from $u_i$ to $v_i$ in $D_i$ with $|L_i| \leq 12 \nu^{-2}$.  and with alternating path $P_i$ satisfying $|P_i| \leq 8\nu^{-1}$. 
By our choice of $D_i$, we see that $L_1, \ldots, L_i$ are vertex-disjoint ladders disjoint from $S_{i+1}$, where $L_i$ is from $u_i$ to $v_i$, completing the induction step and the proof.
\end{proof}

Finally we combine our various constructions to show how to build a $d$-absorber in a robustly expanding digraph.
\begin{theorem}
\label{th:absorber}
Let $0 < \nu \le \tau \le \gamma/16 < 1/16$ and $n,d \in \mathbb{N}$. 
Suppose $d \geq 8$ and
\[
 n>\max\left( 10^4 d^2\gamma^{-5}\log^2(100 d \gamma^{-2}) , 10^5 d \gamma^{-2} \nu^{-3} \log( 1500 d\gamma^{-2} \nu^{-1}) \right).
\]
If $D$ is a robust $(\nu, \tau)$-expander on $n$ vertices with $\delta^{0}(D) \geq \gamma n$ then we can find a $d$-absorber $S$ in $D$ such that $|V(S)| \leq 1600 \nu^{-2} \gamma^{-2} ( d \log(d\gamma^{-2}) + \log n)$.
\end{theorem}
\begin{proof}
For our choice of $\gamma,d,n$, we can apply Lemma~\ref{le:cover} to $D$ to obtain a vertex-disjoint $K \subseteq V(D)^{[2]}$ which $d$-covers $V(D)$, and moreover $m:= |K| = \lceil 24 \gamma^{-2} ( d \log( 24 d\gamma^{-2}) + 2\log n) \rceil$.

Next, by our choice of $n$, we can apply Lemma~\ref{le:ladderscon}
\endnote{
Let $c' = 24000 \gamma^{-2} \nu^{-3}$ and $d' =  12000 \gamma^{-2} \nu^{-3} d \log( 24 d\gamma^{-2})$.
Since $d \ge 8$, we have 
\begin{align*}
		n/2 & \ge 5 \cdot 10^4  \gamma^{-2} \nu^{-3} d \log( 1500 d \gamma^{-2} \nu^{-1}) 
		\ge  4 \cdot 10^5 \gamma^{-2} \nu^{-3}  \log ( 1500 \gamma^{-1/2} \nu^{-3/4} )\\
		& \ge  10^5  \gamma^{-2} \nu^{-3} \log ( 24000 \gamma^{-2} \nu^{-3} )
		 \ge 4 c' \log c'  \ge 3 c' (\log c' +1).
\end{align*}
Also, we have 
\begin{align*}
	 n/2 & \ge 5 \cdot 10^4  \gamma^{-2} \nu^{-3} d \log( 1500 d \gamma^{-2} \nu^{-1})  \ge  36000 \gamma^{-2} \nu^{-3}  \log( 24 d\gamma^{-2}) = 3d'.
\end{align*}
Proposition~\ref{pr:polyexp} implies that
\begin{align*}
	n &  \ge c' \log n +d'
	\ge 12000 \gamma^{-2} \nu^{-3} (  2\log n + d \log( 24 d\gamma^{-2})) \\
	& \ge 460\nu^{-3} \cdot 25 \gamma^{-2}  (  2\log n + d \log( 24 d\gamma^{-2})) 
	\ge 460 m \nu^{-3}.
\end{align*}
}
(taking $k=m$) to construct a ladder from $a$ to $b$ for every $(a,b) \in K$ such that each ladder has at most $12 \nu^{-2}$ vertices, the alternating path of each ladder has length at most $8 \nu^{-1}$, and the ladders are vertex-disjoint. 

Let $\mathcal{L} = \{ L_1, \ldots, L_m \}$ be the set of constructed ladders and let $R_1, \ldots, R_s$ be the collection of all rung paths of all the ladders constructed; thus $s \leq  4 \nu^{-1} m$. Let $x_i$ and $y_i$ be the initial and final vertices of $R_i$ and let $D'$ be the digraph obtained from $D$ by deleting all internal vertices of $R_1, \ldots, R_s$. 
So we have deleted at most $12 \nu^{-2} m \le \nu n/2 $ vertices\endnote{We need $n \ge 24 \nu^{-3}$, which is true by the previous note.}.
Then $D'$ is a $(\frac{1}{2}\nu, \frac{32}{31}\tau)$-expander by Proposition~\ref{pr:del}.
By our choice%
	\endnote{Note that $11 \le s \le 4 \nu^{-1} m  $ and $n > 460m \nu^{-3}$ by Note~6. 
Thus, $|D'| \ge n - 12 \nu^{-2} m \ge 28 \cdot 4 \nu^{-3} m \ge 4 (6 s +11) \nu^{-2}$.}
of parameters and $n$, we can apply Lemma~\ref{le:paths} (with $n=|D'|$, $r=s$ and $ \nu, \tau, \gamma$ replaced by $ \frac{1}{2}\nu, \frac{32}{31}\tau, \frac{31}{32}\gamma$) to find paths $U_i$ from $y_i$ to $x_{i+1}$ for each $i=1, \ldots, s$, where indices are understood to be modulo $s$ and each path has length at most $4\nu^{-1} + 1$. Then $C=x_1R_1U_1 \cdots R_sU_sx_1$ is a cycle in which all the ladders $L_1, \ldots, L_m$ are embedded.
Thus $S = (K, \mathcal{L}, C)$ is a $d$-absorber of $D$.
Also $|V(S)| \leq 12 \nu^{-2}m + s(4\nu^{-1}+1) \leq 32 \nu^{-2} m$.
Recall that $m = \lceil 24 \gamma^{-2} ( d \log( 24 d\gamma^{-2}) + 2\log n) \rceil \le 25 d\gamma^{-2} \log( 24 d\gamma^{-2}) + 48 \gamma^{-2}\log n $ as $d\gamma^{-2} \log( 24 d\gamma^{-2})>1$.
Therefore $|V(S)| \leq 1600 \nu^{-2} \gamma^{-2} ( d \log(d\gamma^{-2}) + \log n)$ as required. 
\end{proof}

\section{Rotation-extension: $1$-factors with few cycles\label{sec:r-e}}

Let $D$ be a digraph. Throughout this section, a \emph{factor} $U$ of $D$ refers to a \emph{$1$-factor} of $D$, i.e.\ a spanning subgraph of $D$ in which every vertex has in- and outdegree $1$. Thus a factor consists of a collection of vertex-disjoint cycles. We shall think of $U$ interchangeably as both a set of vertex-disjoint cycles $U=\{C_1, \ldots, C_k\}$ and as the corresponding subgraph $U= C_1 \cup \cdots \cup C_k$ of $D$. The purpose of this section is to show that any robustly expanding digraph with sufficiently high minimum in- and outdegree contains a factor with few cycles: our main tool is an interesting variation of the rotation-extension technique of P\'{o}sa~\cite{Posa}. The first lemma shows that any robustly expanding digraph with large enough minimum in- and outdegree has a factor.

\begin{lemma}
\label{le:factor}
Let $0 <  \nu \le  \tau < \gamma <1$ and $n \in \mathbb{N}$. If $D=(V,E)$ is an $n$-vertex robust $( \nu, \tau)$-expander with $\delta^{0}(D) \ge \gamma n$ then $D$ has a $1$-factor.
\end{lemma}
\begin{proof}
Let $V=\{v_1, \ldots, v_n\}$. Consider the bipartite (undirected) graph $G$ whose vertex set is $X \cup Y$ where $X = \{x_1,  \ldots, x_n \}$ and $Y = \{y_1, \ldots, y_n \}$ and $x_iy_j$ is an edge of $G$ if and only if $v_iv_j \in E$. Note that $D$ contains a factor if and only if $G$ has a perfect matching, so it is sufficient for us to verify Hall's condition for $G$.
Indeed suppose $S \subseteq X$. If $|S| \leq \tau n$, then $|N_G(S)| = |N_D^+(S)| \geq \gamma n > \tau n \geq |S|$. If $|S| \geq (1- \tau)n$ then since every vertex in $Y$ has degree at least $\gamma n > \tau n$ (since $\delta^-(D) \ge \gamma n$) then $|N_G(S)| = |Y| \ge |S|$.
 If $\tau n \leq |S| \leq (1-\tau) n$, then $|N_G(S)| = |N_D^+(S)| \geq |\RN^+_{\nu}(S)| \geq |S| + \nu n > |S|$. Hence by Hall's Theorem (see e.g. \cite{BondyMurty}) $G$ has a perfect matching and hence $D$ has a factor. 
\end{proof}

We now introduce various notions we shall need.
We say $F$ is a \emph{prefactor} of $D$ if $F$ can be obtained from a factor of $D$ by deleting one edge. Thus $F'$ consists of a collection of cycles $C_1, \ldots, C_{k-1}$ together with a path $P$. We interchangeably think of $F$ as the set $F=\{C_1, \ldots, C_{k-1}, P\}$ and as the subgraph $F = C_1 \cup \cdots \cup  C_{k-1} \cup P$ of $D$. If $P$ is a path from a vertex $x$ to a vertex $y$, we say  $x$ is the \emph{origin} of $F$ and $y$ is the \emph{terminus} of $F$ written $x = {\rm ori}(F)$ and $y = {\rm ter}(F)$ respectively. Every vertex $v$ of $D$ except ${\rm ori}(F)$ has a unique inneighbour in $F$ which we denote by $F^-(v)$.

An \emph{extension} of $F$ (in $D$) is a prefactor $F'$ of $D$ obtained from $F$ as follows. Assuming $F=\{C_1, \ldots, C_{k-1},P\}$, $x = {\rm ori}(F)$ and $y = {\rm ter}(F)$, we pick any vertex $z \in N^+_D(y) \setminus \{ x \}$:
\begin{itemize}
\item[(i)] if $z \in V(P)$ we set $F' = \{C_1, \ldots, C_{k-1}, C',P'\}$, where $C' = zPyz$ and $P' = xPz^-$, where $z^- := F^-(z)$ is the predecessor of $z$ on $P$;
\item[(ii)] if $z \in V(C_i)$ for some $i$ then set $F' = \{C_1, \ldots, C_{i-1}, C_{i+1}, \ldots, C_{k-1}, P' \}$, where $P'=xPyzC_iz^-$ and $z^-:=F^-(z)$ is the predecessor of $z$ in $C_i$.
\end{itemize}
We say $F'$ is an extension of $F$ along the edge $yz$.
Notice that $F$ and $F'$ differ only in their path and in that one or the other contains an additional cycle.
For case (i), we say $F'$ is a \emph{cycle-creating} extension of $F$ and for case (ii) we say $F'$ is a \emph{cycle-destroying} extension of $F$. Notice also that for any extension $F'$ of $F$, we have ${\rm ori}(F) = {\rm ori}(F')$ and that $F'$ is uniquely determined from $F$ by specifying the terminus of $F'$.

Here is the main step in obtaining a factor with few cycles.

\begin{lemma}
\label{le:prefactor}
Let $n \in \mathbb{N}$ and $ \nu, \tau, \gamma, \xi \in (0,1)$ satisfying $ \nu \le \tau$, $\gamma > 2\tau + \xi$, $\xi < \frac{1}{4}\nu^2$ and $n > 32 \nu^{-3}$.
Suppose $D=(V,E)$ is an $n$-vertex robust $(\nu, \tau)$-expander with $\delta^{0}(D) \geq \gamma n$ and suppose that for each prefactor $F$ of $D$, we have an associated set $B(F) \subseteq V$ of `forbidden' vertices satisfying ${\rm ori}(F) \in B(F)$ and $|B(F)| \leq \xi n$.
Fix any prefactor $F^*$ of $D$. Then for all but at most $\tau n$ vertices $y \in V$, there exists a sequence of prefactors $F_0=F^*, F_1, \ldots, F_t$ where $y = {\rm ter}(F_t)$ and for each $i = 1, \ldots, t$ we have that $F_i$ is an extension of $F_{i-1}$ and ${\rm ter}(F_i) \not\in B(F_{i-1})$.
\end{lemma}
\begin{proof}
Let $x = {\rm ori}(F_0) = {\rm ori}(F^*)$. For each $r \in \mathbb{N}$, we define $S_r$ to be the set of vertices that are \emph{reachable} from $F_0$ by a sequence of at most $r$ successive extensions while avoiding forbidden sets. More precisely, $y \in S_r$ if and only if there exists a sequence $F^*=F_0, F_1, \ldots, F_{r'}$ with $r' \leq r$ such that $y = {\rm ter}(F_{r'})$, and for all $i = 1, \ldots, r'$, $F_i$ is an extension of $F_{i-1}$ and ${\rm ter}(F_i) \not\in  B(F_{i-1})$. For each $y \in S_r$, we set $F^{(r)}_y := F_{r'}$
(if there are many choices of $F_{r'}$, we pick one arbitrarily). In particular $y \in {\rm ter}(F^{(r)}_y)$.

In order to prove the lemma, it is sufficient to show that $|S_t| \geq (1 - \tau)n$ for some $t$. Let us begin by noting that $|S_1| \geq (\gamma - \xi)n - 1 \geq 2\tau n - 1 \geq \tau n$, where the last two inequalities follow by our choice of parameters and $n$. To see the first inequality note that each distinct outneighbour $w$ of ${\rm ter}(F_0)$ (except possibly $x$) gives an extension of $F_0$ with a distinct terminus $w^-:=F_0^-(w)$, and each such $w^-$ is in $S_1$ unless $w^- \in B(F_0)$.

We shall show that $S_{r+1}$ contains most vertices in $\{ F_0^-(w): w \in {\rm RN}_{\nu}^+(S_r)\}$.
Fix $r \geq 1$. For each $w \in {\rm RN}^+_{\nu}(S_r) \setminus \{ x \}$, we say $w$ is \emph{good} if there exists $v \in S_r$ such that $w \in N^+(v)$, $F_v^{(r)-}(w) = F_0^-(w)$, and $F_0^-(w) \not\in B(F^{(r)}_v)$. Otherwise we say $w$ is \emph{bad}. Note that if $w \in {\rm RN}^+_{\nu}(S_r) \setminus \{ x \}$ is good, then $F_0^-(w) \in S_{r+1}$. Indeed, let $F_0, \ldots, F_{r'}=F^{(r)}_v$ be a sequence of extensions that show $v \in S_r$. Then extending $F_{r'}=F^{(r)}_v$ along the edge $vw$ gives an extension $F'$ whose terminus is $F_v^{(r)-}(w) = F_0^-(w) \not\in B(F^{(r)}_v)$. Thus the sequence $F_0, \ldots, F_{r'}, F'$ shows that $F_0^-(w) \in S_{r+1}$.

Since the function $w \mapsto F_0^-(w)$ is injective, each $w \in {\rm RN}^+_{\nu}(S_r) \setminus \{ x \}$ that is good corresponds to a distinct vertex of $S_{r+1}$. Thus, assuming $|S_r| \leq (1- \tau) n$, we have
\[
|S_{r+1}| \geq |{\rm RN}^+_{\nu}(S_r)| - b - 1 \geq |S_r| + \nu n -b - 1,
\]
where $b$ is the number of bad vertices, which we now bound from above.

Let 
\begin{align*}
A &:= \{(v,w):\; v \in S_r, \; w \in \RN^+_{\nu}(S_r) \cap N^+(v),\; F_0^-(w) \in B(F^{(r)}_v) \} \\
B &:= \{(v,w):\; v \in S_r, \; w \in \RN^+_{\nu}(S_r) \cap N^+(v),\; F^{(r)-}_v(w) \not= F_0^-(w) \}.
\end{align*}
We have that $|A \cup B| \geq b \nu n$. To see this note that 
each bad vertex $w \in \RN^+_{\nu}(S_r)$ has at least $\nu n$ inneighbours $v \in S_r$, and each such pair $(v,w)$ belongs to $A \cup B$. On the other hand, we have $|A| \leq \sum_{v \in S_r}|B(F^{(r)}_v)| \leq  |S_r| \xi n$ and $|B| \leq |S_r|r$. The first inequality is clear while second inequality follows from the following claim:

\medskip
\noindent
{\bf Claim}: For each $v \in S_r$, there are at most $r$ vertices $w$ for which $F^{(r)-}_v(w) \not= F_0^-(w)$.
\begin{proof}(of Claim)
If $F'$ is any extension of $F$ then there is exactly one vertex $w$ for which $F'^-(w) \not= F^-(w)$.
Therefore if $F'$ is obtained from $F$ by a sequence of at most $r$ successive extensions, then there are at most $r$ vertices $w$ for which $F'^-(w) \not= F^-(w)$, and the claim follows.
\end{proof}  

Thus we have that $b \nu n \leq (\xi n + r) |S_r|$, whence $b \leq \nu^{-1} \xi n + \nu^{-1}r$. For each $r \leq 2 \nu^{-1}$, if $|S_r| \leq (1 - \tau)n$ then we have
\[
|S_{r+1}| \geq |S_r| + \nu n - b -1 \geq |S_r| + \nu n - \nu^{-1} \xi n - 2 \nu^{-2} - 1 \geq |S_r| + \frac{1}{2}\nu n,
\]
where the last inequality follows\endnote{We note  $\nu^{-1} \xi n \leq \frac{1}{4} \nu n$ and $1 \leq 2 \nu^{-2} \leq \frac{1}{8} \nu n$.} by our choice of parameters and $n$. Thus for some $t \leq 2 \nu^{-1}$, we have $|S_t| \geq (1 - \tau)n$, as required.
\end{proof}

We give one piece of notation before proving the existence of factors with few cycles in robustly expanding digraphs. If $P$ and $Q$ are paths in a directed graph $D$, we write $Q \subseteq P$ if $Q$ is an initial segment of $P$, i.e.\ $P$ and $Q$ have the same initial vertex and $P[V(Q)] = Q$. If $Q \subseteq P$ but $Q \not= P$, we write $Q \subset P$.

\begin{theorem}
\label{th:factor}
Let $n \in \mathbb{N}$ and $ \nu, \tau, \gamma, \xi \in (0,1)$ satisfying $ \nu \le \tau$, $\gamma > 2\tau + \xi$, $\xi < \frac{1}{4}\nu^{2}$, and $n > 32 \nu^{-3}$.
If $D=(V,E)$ is an $n$-vertex robust $(\nu, \tau)$-expander with $\delta^{0}(D) \geq \gamma n$ then there exists a factor $U^*$ of $D$ which consists of at most $2\xi^{-1}$ cycles.
\end{theorem}
\begin{proof}
By Lemma~\ref{le:factor}, $D$ contains a factor $U_0$. Suppose $U$ is any factor in which all cycles have length at least $s$ for some $s < \frac12 \xi n$ and where exactly $\ell \geq 1$ cycles have length $s$. We claim that, using Lemma~\ref{le:prefactor}, we can obtain a factor $U'$ from $U$ in which all cycles have length at least $s$ and at most $\ell -1$ cycles have length $s$. Applying this claim iteratively, we eventually obtain a factor $U^*$ of $D$ in which every cycle has length at least $\frac{1}{2}\xi n$ and so this factor has at most $2\xi^{-1}$ cycles, proving the theorem.

It remains to prove the claim. Suppose $U = \{C_1, \ldots, C_k \}$ where $C_1, \ldots, C_k$ are the cycles of $U$ in increasing order of length with $|C_1| = s < \frac12 \xi n$. Delete any edge of $C_1$ to form a path $P$ and let $F_0 = \{C_2, \ldots, C_k, P \}$ be the resulting prefactor of $D$, and let $x$ be its origin.

For each prefactor $F$ of $D$, let $B(F)$ denote the set of the first and last $\frac{1}{2}\xi n$ vertices on the path in $F$ (if the path has at most $\xi n$ vertices then $B(F)$ is the set of all vertices on the path). Note that for the prefactor $F_0$, $|P| = |C_1| < \frac12 \xi n$ and so $B(F_0) = V(P)$. 
By Lemma~\ref{le:prefactor}, for at least $(1- \tau)n$ vertices $y \in V$, there exists a sequence of extensions $F_0, F_1, \ldots, F_t$ such that $F_i$ is an extension of $F_{i-1}$, ${\rm ter}(F_i) \not\in B(F_{i-1})$, and ${\rm ter}(F_t) = y$. Since $|N^-(x) \setminus B(F_0)| \geq \gamma n - \xi n > \tau n$, we can choose $y$ to be in $N^-(x) \setminus B(F_0)$.

Writing $P_i$ for the path in the prefactor $F_i$, by our choice of $B(\cdot)$, it is straightforward to show by induction that $P = P_0 \subset P_i$ for all $i = 1, \ldots, t$. Indeed, since $B(F_0) = V(P)$, $F_1$ must be a cycle-destroying extension of $F_0$, and so $P = P_0 \subset P_1$. Suppose $P \subset P_{i-1}$ for some $i>1$ and let $P_{i-1}'$ be the subpath of $P_{i-1}$ consisting of the first $\frac{1}{2}\xi n$ vertices; in particular $P \subset P_{i-1}'$. If $F_i$ is a cycle-creating extension of $F_{i-1}$, then since $V(P_{i-1}') \subseteq B(F_{i-1})$, we must have $P_i \supseteq P_{i-1}' \supset P$. If $F$ is a cycle-destroying extension of $F_i$ the $P_i \supset P_{i-1} \supset P_0$.

Our choice of $B(\cdot)$ also ensures that if $F_i$ is a cycle-creating extension of $F_{i-1}$, then the new cycle has length at least $\frac{1}{2}\xi n$.

Let $F_t = \{C_1' , \ldots, C_{k'}', P_t \}$, where $C_1' \ldots, C_{k'}'$ are cycles and we know $P_t$ is a path from $x$ to $y$ of length more than $|P|= |C_1|$. Since $y \in N_D^-(x)$, we can turn $P_t$ into a cycle $C^*$ and form a factor $U' = \{C_1' \ldots, C_{k'}', C^* \}$ of $D$. We have $|C^*| = |P_t| > |C_1| = s$.

 Every cycle in $U'$ that was created in the sequence of extensions $F_0, \ldots, F_t$ has length at least $\frac12 \xi n > s$ and $|C^*| > |C_1|=s$. Every other cycle of $U'$ was also a cycle of $U$. Hence every cycle in $U'$ has length at least $s$ and the number of cycles of length exactly $s$ has been reduced by at least one. This proves the claim and the theorem. 
\end{proof}

\section{Hamiltonicity}
\label{sec:Ham}

We now combine Theorem~\ref{th:absorber}, Corollary~\ref{co:d-absorb} and Theorem~\ref{th:factor} to give the following result from which we deduce Theorem~\ref{th:maino}. 

\begin{theorem}
\label{th:main}
Let $0 <  \nu \le \tau \le \gamma/16 <1/16$ and let $n \in \mathbb{N}$. Assume 
\[
n> \max \{  10^8 \gamma^{-5} \nu^{-4} \log^2(10^4 \gamma^{-2}\nu^{-2}) , 10^7 \gamma^{-2} \nu^{-5} \log( 150000 \gamma^{-2} \nu^{-3}) \}.
\]
If $D$ is an $n$-vertex robust $(\nu,\tau)$-expander with $\delta^{0}(D) \geq \gamma n$, then for any $\nu  n/2 \le  \ell \le n $ and any $v \in V(D)$, $D$ contains a cycle of length $\ell$ containing~$v$.
\end{theorem}
\begin{proof}
Let $\xi := \nu^2 /32$ and $d:= \lceil 2 \xi^{-1} \rceil \ge 8$.
We begin by applying Theorem~\ref{th:absorber} to $D$ to find a $d$-absorber $S$, where 
\begin{align*}
|V(S)| &\leq  1600 \nu^{-2} \gamma^{-2} ( d \log(d\gamma^{-2}) + \log n).
\end{align*}
One can check that the conditions on the parameters and $n$ are met%
\endnote{
Note that $d := \lceil 2 \xi^{-1} \rceil \leq 3 \xi^{-1} \leq 100 \nu^{-2}$ and $\nu^{-1} \ge 16$.
Hence
\begin{align*}
10^4 d^2\gamma^{-5} \log^2(100 d\gamma^{-2}) \le 10^8 \gamma^{-5} \nu^{-4} \log^2(10^4 \gamma^{-2}\nu^{-2}) \le n 
\end{align*}
and 
\begin{align*}
10^5 d \gamma^{-2} \nu^{-3} \log( 1500 d\gamma^{-2} \nu^{-1})
\le  10^7 \gamma^{-2} \nu^{-5} \log( 150000 \gamma^{-2} \nu^{-3}) \le n .
\end{align*}
}.

Set $D' := D - V(S)$. By our choice%
	\endnote{
	Need $n > 2\nu^{-1}|V(S)|$ so sufficient that $n > 3200 \gamma^{-2} \nu^{-3} ( d \log(d\gamma^{-2}) + \log n)$.
	By Proposition~\ref{pr:polyexp} this holds if 
\[
n > 9600 \gamma^{-2} \nu^{-3} \left( \log (3200 \gamma^{-2} \nu^{-3})+1\right) + 9600  \gamma^{-2} \nu^{-3} d \log(d\gamma^{-2}).
\] 
Recall that $d \leq 100 \nu^{-2}$.
The inequality above holds if $n > 10^7  \gamma^{-2}  \nu^{-5} \log(10^5   \gamma^{-2} \nu^{-4})$, which holds if $n > 3 \cdot 10^7  \gamma^{-2}  \nu^{-5} \log(   \gamma^{-2} \nu^{-4} )$.
} 
of $n$, we have $|V(S)| < \nu n/2$ and so by Proposition~\ref{pr:del} $D'$ is a robust $( \frac{1}{2}\nu, \frac{32}{31}\tau)$-expander with $\delta^{0}(D') > \frac{31}{32}\gamma n$. 
By our choice of $\xi$ and $n$, we can apply Theorem~\ref{th:factor}%
	\endnote{We check that $\frac{31}{32}\gamma > \frac{64}{31}\tau + \xi$, which holds (using $\gamma > 16\tau$ and $\xi \le \nu \le \tau$).
We check that $\xi < \frac{1}{4}(\frac{1}{2}\nu)^2 = \frac{1}{16} \nu^{2}$, which holds. 
We check $n - |V(S)| > \max(32(\frac{1}{2}\nu)^{-3}, \tau^{-1})$. 
Since $|V(S)| \leq \frac{1}{2} \nu n \leq \frac{1}{2}n$, it is sufficient that $n >  512\nu^{-3}$. This is clearly implied by our choice of $n$.
} 
to $D'$ to obtain a factor in $D'$ with at most $2 \xi^{-1} \le d$ cycles.
By removing one edge from each of the cycles let $P_1, \ldots, P_r$ be the resulting paths with $r \le d$.
Consider any $\nu n/2 \le \ell \le n$ and any $v \in V(D)$. 
Note $D'$ contains vertex-disjoint paths $P'_1, \ldots, P'_{r'}$ such that $r' \le d$ and $|P'_1| + \ldots + |P'_{r'}| = \ell - |V(S)|$ and $v \in V(S) \cup \bigcup_{i \in r'} V(P'_i) $ (by removing appropriate vertices of $P_1, \ldots, P_r$ if necessary).
Applying Corollary~\ref{co:d-absorb}, to these paths and the $d$-absorber $S$, we obtain a cycle $C$ of length $\ell$ in~$D$ with $v \in V(C)$.
\end{proof}

Finally we can prove Theorem~\ref{th:maino}.

\begin{proof}[Proof of Theorem~\ref{th:maino}]
Given that $D$ is an $n$-vertex robust $(\nu, \tau)$-outexpander with $\delta^0(D) \geq \gamma n$, by Proposition~\ref{pr:dib}, $D$ is a robust $(\nu',\tau)$-expander where $\nu' = \nu^2/2$ (our choice of parameters ensures the conditions of Proposition~\ref{pr:dib} are met).
By our choice
	\endnote{
	Recall that $\nu \le \gamma/16$. so
	\begin{align*}
		 & \max \{  10^8 \gamma^{-5} (\nu')^{-4} \log^2(10^4 \gamma^{-2} (\nu')^{-2}) , 10^7 \gamma^{-2} (\nu')^{-5} \log( 15000 \gamma^{-2} (\nu')^{-3}) \}\\
		& \le 40000 \nu^{-13}		\log^2 250 \nu^{-8} \le 40000 \nu^{-13}		\log^2 \nu^{-10} \le   (4 \nu^{-1})^{13}		\log^2 \nu^{-1} 
\le n .
\end{align*}
}
 of $n$ we can apply Theorem~\ref{th:main} to $D$ to obtain a Hamilton cycle in~$D$. 
\end{proof}

We deduce Corollary~\ref{co:luke} from Theorem~\ref{th:maino}, but first we need the following leema from~\cite[Lemma~13.1]{Kelly}. 

\begin{lemma}
\label{le:orient}
Let $n \in \mathbb{N}$ and $\nu, \tau, \varepsilon \in (0,1)$ satisfy $\nu \le \frac{1}{8}\tau^2$ and $\tau \le \frac{1}{2} \varepsilon$. If $D$ is an oriented graph on $n$ vertices with $\delta^+(D) + \delta^-(D) + \delta(D) \geq 3n/2 + \varepsilon n$ then $G$ is a robust $(\nu, \tau)$-outexpander. 
\end{lemma} 
The explicit dependence between the parameters was not given in \cite{Kelly}, but we have computed them and included them in the statement above.

\begin{proof}[Proof of Corollary~\ref{co:luke}]
Let $\nu = \varepsilon^2 / 2$, $\tau := 2 \varepsilon$ and $\gamma = 3/8$ so $4 \sqrt[13]{\log^2 n / n } \le \nu \le \tau \le \gamma /16 <1/16$.
Given an $n$-vertex oriented graph $D$ with $\delta^0(D) \geq 3n/8 + \varepsilon n$, we have that $\delta(D) + \delta^+(D) + \delta^-(D) \geq 3n/2 + 4 \varepsilon n$.
So by Lemma~\ref{le:orient},  $D$ is a robust $(\nu, \tau)$-outexpander.
Finally, we apply Theorem~\ref{th:maino} to obtain a Hamilton cycle.
 
\end{proof}

To prove Theorem~\ref{NWa}, we need the following lemma from~\cite[Lemma~13]{KOT2}.

\begin{lemma}
\label{le:orient2}
Let $n \in \mathbb{N}$ and $\tau, \gamma \in (0,1)$ satisfy $2\tau + 4 \tau^2 \le \gamma \le 1/2$ and $n \ge \gamma^{-2}$.
Let $D$ be an $n$-vertex digraph such that for all $i < n/2$,
\begin{itemize}
	\item $d_i^+ \ge i + \gamma n$ or $d_{n-i-\gamma n }^- \ge n-i$,
	\item $d_i^- \ge i + \gamma n$ or $d_{n-i-\gamma n }^+ \ge n-i$.
\end{itemize}
Then $G$ is a robust $(\tau^2, \tau)$-outexpander and $\delta^0(D) \ge \gamma n$.
\end{lemma} 

\begin{proof}[Proof of Theorem~\ref{NWa}]
Let $\tau := \gamma/16$.
By~Lemma~\ref{le:orient2}, $D$ is a robust $(\tau^2, \tau)$-outexpander and $\delta^0(D) \ge \gamma n$. Finally, we apply Theorem~\ref{th:maino} to obtain a Hamilton cycle.
\end{proof}

\section{Concluding remarks and an open problem}
\label{sec:conc}

It would be interesting to know for which choices of parameters $\nu = \nu(n)$ and  $\tau = \tau(n)$ an $n$-vertex robust $(\nu, \tau)$-expander is guaranteed to be Hamiltonian. We believe the true values of $\nu$ and $\tau$ for which this holds should be much smaller than what we have proved.

\subsection*{Acknowledgements}
The authors would like to thank the anonymous referee for their careful reading and helpful suggestions.

\theendnotes

\end{document}